\documentclass[10pt]{article}
\usepackage{amsfonts}
\usepackage{}
\usepackage{bm}
\usepackage{amsthm}

\usepackage{amsmath}
\usepackage{amssymb}
\usepackage{graphicx}
\usepackage{epic,eepic}

\usepackage{color}

\usepackage{amstext,amsthm,amssymb,amsmath}
\usepackage{graphicx}
\usepackage{subfigure}
\usepackage{wrapfig}
\usepackage{fullpage}
\usepackage{color}
\usepackage{multirow}
\usepackage{tabulary}
\usepackage{booktabs}
\usepackage{enumerate}
\usepackage{epsfig,epstopdf}

\usepackage{soul}

\newcommand{\om}{\Omega}

\newtheorem{theorem}{Theorem}[section]

\newtheorem{remark}[theorem]{Remark}

\newcommand{\eps}{\varepsilon}

\newcommand{\avrg}[2]{{\left \langle #1 \right \rangle}_{#2}}

\hfuzz=\maxdimen
\graphicspath{ {./pic/}}

\title{Generalized Multiscale Finite Element Methods for problems in  perforated heterogeneous domains}

\author{
Eric T. Chung \thanks{Department of Mathematics,
The Chinese University of Hong Kong (CUHK), Hong Kong SAR. Email: {\tt tschung@math.cuhk.edu.hk}.
The research of Eric Chung is supported by Hong Kong RGC General Research Fund (Project 400411).}
\and
Yalchin Efendiev \thanks{Department of Mathematics \& Institute for Scientific Computation (ISC),
Texas A\&M University,
College Station, Texas, USA
and Center for Numerical Porous Media (NumPor),
King Abdullah University of Science and Technology (KAUST),
Thuwal 23955-6900, Kingdom of Saudi Arabia. Email: {\tt efendiev@math.tamu.edu}.}
\and
Guanglian Li\thanks{Department of Mathematics, Texas A\&M University, College Station, TX 77843-3368} \and
 Maria Vasilyeva\thanks{Department of Computational Technologies, Institute of Mathematics and Informatics, North-Eastern Federal University, Yakutsk, Republic of Sakha (Yakutia), Russia, 677980 \& Institute for Scientific Computation, Texas A\&M University, College Station, TX 77843-3368}
}
\begin{document}
\maketitle
%==============================

\begin{center}
{\it Submitted to the Special Issue Mathematical \& Numerical Analysis of Flow and Transport in Porous Media}
\end{center}

\begin{abstract}

Complex processes in perforated domains occur in many real-world applications.
These problems are typically characterized by physical processes in domains with multiple
scales (see Figure \ref{fig:perf_domain} for the illustration of a perforated domain). Moreover,
these problems are intrinsically multiscale and their discretizations can yield very large
linear or nonlinear systems. In this paper, we investigate multiscale approaches that attempt
to solve such problems on a coarse grid by constructing multiscale basis functions in each coarse
grid, where the coarse grid can contain many perforations. In particular, we are interested in cases
when there is no scale separation and the perforations can have different sizes.
In this regard, we mention some earlier pioneering works \cite{Henning09,Bris14, CR13}, where the authors
develop multiscale finite element methods. In our paper, we follow Generalized Multiscale Finite Element
Method (GMsFEM) and develop a multiscale procedure where we identify multiscale basis functions
in each coarse block using snapshot space and local spectral problems. We show that with a few
basis functions in each coarse block, one can accurately approximate the solution, where each coarse
block can contain many small inclusions. We apply our general concept to
(1) Laplace equation in perforated domain; (2) elasticity equation in
perforated domain; and (3) Stokes equations in perforated domain.
Numerical results are presented
for these problems using two types of heterogeneous perforated domains.
The analysis of the proposed methods will be presented elsewhere.

\end{abstract}

\section{Introduction}

%Conclusions.

%motivation for perforated domain.
Among multiscale problems, the problems in perforated domains
are of great interest for many applications.
The main characteristics of these problems is that the underlying processes
occur in multiscale domains where the geometry of the domain has multiple scales,
e.g., domain outside inclusions (see Figure \ref{fig:perf_domain}).
%Partial differential operators in perforated domain is of great interest in the past
%few years.
There are many important applications for processes in perforated domains.
For example, fluid flow in porous media, diffusion in perforated domains,
mechanical processes in hollow materials, and so on.
Typically, it is a combination of a physical process and a heterogeneous
media that gives rise to problems in perforated domains, e.g.,
fluid flow in porous media is a problem in perforated domain, while heat conduction
in porous media may not be a problem in perforated domain.

%General multiscale and the problems in perforated domains.
The problems in perforated domains are of multiscale nature.
%require large degrees of freedom to discretize and solve them.
The solution techniques for these problems
require high resolution. In particular, the discretization needs
to honor the irregular boundaries of perforations. This gives rise
to a fine-scale problems with many degrees of freedom which can
be very expensive to solve. In this paper, our goal is to develop
coarse-grid approaches where the coarse grids do not have to align
with perforations and the perforated domains do not need to have a scale
separation. Our objective is to construct a low dimensional approximate model
for solving fine-scale problems.

%Challenges in developing multiscale methods in perforated domains and previous work on multiscale.
When the perforated domains have some scale separation,
there are numerous works which include the works on the homogenization and
asymptotic expansion in
periodic perforated domains \cite{Jikov91, Cao06, Cao10}.
Recently, a novel method of meso-scale asymptotic approximation for
Laplace operator is introduced (\cite{Maz'ya13}) for the cases with a
large number of perforations. The multiscale analysis of the
perforated domain has been conducted using heterogeneous multiscale
finite element method (HMMFEM) and Mulitiscale Finite Element Method (MsFEM)
\cite{Jikov91, Henning09,Bris14}. The authors
in \cite{CR13, Bris14, Muljadi14} have extended
multiscale finite element methods to
arbitrary perforated
domains and presented novel numerical approaches using Crouzeix-Raviart coupling of
 multiscale finite
element basis. Their approach allows perforations to intersect the boundaries
of coarse-grid blocks.
In their approach, the authors use oversampling technique and weak coupling via Crouzeix-Raviart
to avoid directly imposing boundary conditions for the basis functions. The latter
can be difficult for multiscale basis function construction in perforated domains.
In our paper,  we also follow this general concept and avoid constructing boundary conditions
for multiscale basis functions and construct them via local spectral decomposition.
Our approach follows a Generalized Mulitiscale Finite Element Method (GMsFEM)
\cite{egh12} in constructing multiscale basis functions which we discuss next.

%Our work. GMsFEM general concept.
GMsFEM is a general multiscale procedure where
multiscale basis functions in each coarse grid are constructed
and is designed for many applications \cite{WaveGMsFEM,ElasticGMsFEM,MixedGMsFEM}.
In the GMsFEM framework,  one divides the computations into two stages, i.e.,
the offline stage and the online stage. In the offline stage, a reduced dimensional space is constructed,
and it is then used in the online stage to construct multiscale basis functions. These multiscale basis
functions can be re-used for any input parameter to solve the problem on a coarse grid. The main
idea behind the construction of the offline and online spaces is to design appropriate
snapshot spaces and determine an appropriate
local
spectral problem to select important modes in the snapshot space.
In \cite{egh12},
several general strategies for designing the local spectral
 procedures have been proposed. To apply GMsFEM for problems in perforated domains,
one needs a novel way to construct the snapshot space and local spectral problem.
In earlier works \cite{Babuska, hw97, ehw99}, heterogeneous
 problems  have been studied.
The main
challenge is to take into account the domain heterogeneities when constructing multiscale basis functions.
Compared to problems in heterogeneous media without perforated regions, we impose
boundary conditions on perforations. We
generate local snapshot space using
harmonic extensions of appropriately chosen a large set of boundary conditions. Local spectral
problems are used to select dominant modes. We also discuss the use of randomized snapshots
(\cite{cegl14}) that allow calculating a few snapshot vectors. The latter is important in reducing
the computational cost.

%Problems we have considered.

The above procedure is formulated for a general problem in perforated media. In this paper,
we apply it to three different problems that arise in many applications.
These include: (1) Laplace equation in perforated domain; (2) elasticity equation in
perforated domain; and (3) Stokes equations in perforated domain. We show that our general concept
can be applied to all these problems.  We consider several numerical examples using two sets
of perforated domains. One domain contains several large inclusions (with the sizes comparable to
the coarse-grid sizes) and the other domain contains small inclusions (with the sizes
much smaller than the size of the domain). We solve the above three PDEs in these two sets of
domains and compare the multiscale solutions with the fine-grid solutions.
 Our preliminary
numerical results demonstrate a fast convergence of the proposed
 GMsFEM approaches. In particular,
we show that one can achieve a good accuracy with a very few degrees of freedom, especially for
Laplace and elasticity problems. For Stokes equations, there is a room for improvement
by enriching the pressure space. We discuss it in our numerical results.

%Numerical results and some findings.

The rest of the paper is organized as follows. In Section \ref{prelim},
we present
preliminaries on the problems in perforated domains and the GMsFEM framework.
The construction of the coarse spaces for the GMsFEM is discussed
in Section \ref{locbasis}.
In Section \ref{sec:numerical}, numerical results for several
representative examples are presented.
Finally, we conclude our paper with some remarks in
Section \ref{sec:conclusion}.

%structure
\section{Preliminaries}
\label{prelim}
%==============================

\begin{figure}[htb]
  \centering
  \includegraphics[width=0.65 \textwidth]{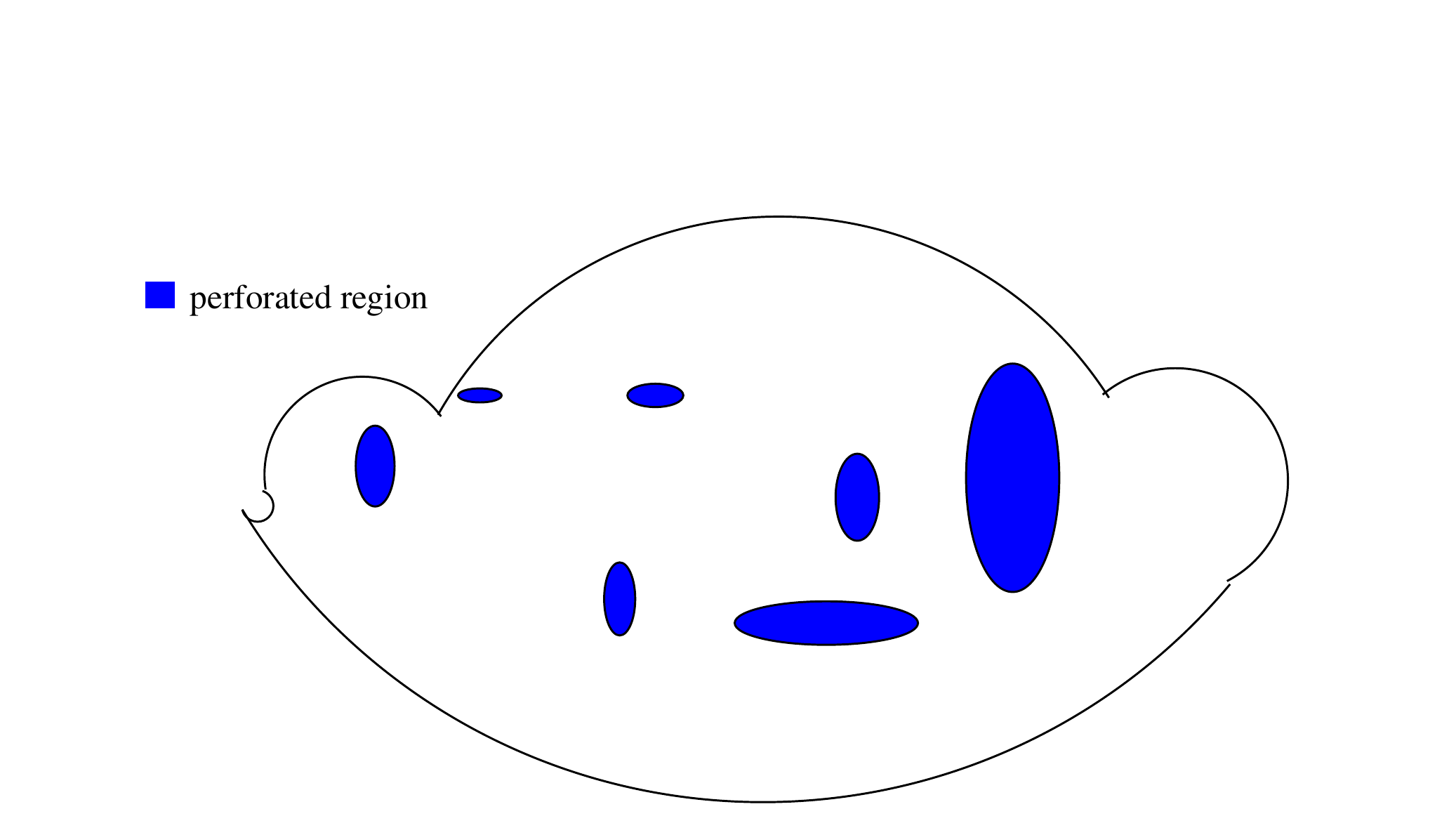}
  \caption{Illustration of a perforated domain.}
  \label{fig:perf_domain}
\end{figure}

In this section, we first present the underlying problem
and the corresponding fine-scale discretization. Then,
we discuss the multiscale strategy of solving this problem.
We consider
\begin{align} \label{eq:original}
&\mathcal {L}^{\eps}(u)=f \quad \text{in} \quad \Omega^{\eps},&\\
&u=0 \text{ or } \frac{\partial u}{\partial n}=0, \text{ on }\partial \Omega^{\eps}\cap  \partial \mathcal {B}^{\eps},&\\
&u=g, \text{ on }\partial \Omega\cap \partial \Omega^{\eps}.&
\end{align}
Here, $\Omega\in \mathbb{R}^d$ ($d=2,3$) is a bounded domain covered
by inactive cells (for Stokes flow and Darcy flow) or active cells (for elasticity problem)
$\mathcal {B}^{\eps}$.
We call active cells where the underlying problem is solved, while inactive cells being the rest of the region.
Suppose the distance between inactive cells (or active cells) is at most $\eps$
and we use $\epsilon$ scripts to denote the perforated domains. Denote the remaining part as $\Omega^{\eps}$, i.e. $\Omega^{\eps}=\Omega\backslash \mathcal {B}^{\eps}$. See Figure \ref{fig:perf_domain} for an illustration of the perforated domain. $\mathcal {L}^{\eps}$ denotes a linear differential operator, e.g. $\mathcal {L}^{\eps}(u)=-\Delta u$ for Darcy flow. $n$ is the unit outward normal to the boundary and $f\in F(\om^{\eps})$ and $g\in G(\om^{\eps})$ denote functions with a suitable regularity.
Note that $\eps$ is used to denote heterogeneous domains, even though
we do not assume periodicity or scale separation.

To simplify the notation, we denote by ${V}(\om^{\eps})$ the appropriate solution space, and \[{V}_{0}(\om^{\eps})=\{v\in {V}(\om^{\eps}), v=0 \text{ on }\partial\om^{\eps}\}.\]
The variational formulation of Problem \eqref{eq:original} is to find $u\in {V}(\om^{\eps})$ such that
\[
\avrg{\mathcal {L}^{\eps}(u),v}{\om^{\eps}}= \avrg{f,v}{\om^{\eps}} \qquad  \text{for all } v \in V_{0}(\om^{\eps}),
\]
$\avrg{\cdot,\cdot}{\om^{\eps}}$ denotes a specific inner product over $\om^{\eps}$ for either scalar functions or vector functions.
%In Example \ref{ex:operators}, different types of differential operators and spaces are defined.
In the following, we give some specific examples for the above abstract notations.
%\begin{example}\label{ex:operators}
\begin{itemize}
\item For the Laplace operator with homogeneous Dirichlet boundary conditions on $\partial \om^{\eps}$, we have
\begin{align}\label{eqn:laplace}
\mathcal {L}^{\eps}(u)=-\Delta u,
\end{align}
and ${V}(\om^{\eps})=H^{1}_{0}(\Omega^{\eps})$, $\avrg{\mathcal {L}^{\eps}(u),v}{\om^{\eps}}=\avrg{\nabla u,\nabla v}{\om^{\eps}}$.

\item For the elasticity operator with homogeneous Dirichlet boundary condition on $\partial \om^{\eps}$,
we let $ {u}\in (H^{1}(\Omega^{\eps}))^{d}$ be the displacement field.
The strain tensor $ {\varepsilon}( {u})\in (L^{2}(\Omega^{\eps}))^{d\times d}$
is defined by
\begin{equation*}
 {\varepsilon }( {u}) = \frac{1}{2} ( \nabla  {u} + \nabla  {u}^T ).
\end{equation*}
In this paper, we assume the medium is isotropic.
Thus, the stress tensor $ {\sigma}( {u})\in (L^{2}(\Omega^{\eps}))^{d\times d}$ relates to the strain tensor $ {\varepsilon}( {u})$ in the following way
\begin{equation*}
 {\sigma}(u)= 2\mu  {\varepsilon} + \lambda \nabla\cdot  {u} \,  {I},
\end{equation*}
where $\lambda>0$ and $\mu>0$ are the Lam\'e coefficients. We have
\begin{align}\label{eqn:elasticity}
\mathcal {L}^{\eps}(u)=- \nabla \cdot  {\sigma},
\end{align}
where ${V}(\om^{\eps})=(H^{1}_{0}(\Omega^{\eps}))^{d}$ and
$\avrg{\mathcal {L}^{\eps}(u),v}{\om^{\eps}}
=\avrg{\sigma(u),\sigma(v)}{\om^{\eps}}$.
\item For Stokes equations, we have
\begin{align}\label{eqn:stokes}
\mathcal {L}^{\eps}(u\;,p)=\begin{pmatrix}
\nabla p -\mu\Delta {u}\\
\nabla \cdot {u}
\end{pmatrix},
\end{align}
where $\mu$ is the viscosity, $p$ is the fluid pressure, $u$ represents the velocity, ${V}(\om^{\eps})=(H^{1}_{0}(\Omega^{\eps}))^{d}\times L^{2}_{0}(\Omega^{\eps})$, and
\[\avrg{\mathcal {L}^{\eps}(u\;,p),(v\;,q)}{\om^{\eps}}=
\begin{pmatrix}
\avrg{\nabla u,\nabla v}{\om^{\eps}} &-\avrg{\nabla \cdot v,p}{\om^{\eps}}\\
\avrg{\nabla \cdot u,q}{\om^{\eps}}&0
\end{pmatrix}.\]
We recall that $L^{2}_{0}(\Omega^{\eps})$ contains functions in $L^{2}(\Omega^{\eps})$
with zero average in $\Omega^{\eps}$.
\end{itemize}
%\end{example}
 We now turn our attention to a numerical
approximation of the variational problem above.
Due to the nature of this multiscale problem we consider the framework of
multiscale finite element method \cite{hw97}. The idea is computing the solution on a coarse grid instead of calculating on the fine grid directly. In the following, we introduce the necessary concepts and notations.

Let $\mathcal{T}_H$ be a coarse-grid partition of the domain
$\Omega^{\eps}$ and $\mathcal{T}_h$ be a
conforming fine triangulation of $\Omega^{\eps}$.
We assume that $\mathcal{T}_h$ is a
refinement of $\mathcal{T}_H$, where $h$ and $H$ represent
the fine and coarse mesh sizes,
respectively. Typically,
 we assume that $0 < h \ll H < 1$, and that the fine-scale mesh $\mathcal{T}_h$
is sufficiently fine to fully resolve the small-scale information of the domain while $H$ is a coarse mesh containing many fine-scale features.
On the triangulation $\mathcal{T}_h$,
we introduce the following finite element spaces %{\blue{for simplex do we still need Q2/Q0? }}
\begin{align*}
{V}_h  &:= \{ {v} \in V(\Omega^{\eps})| {v}|_K \in (P^k(K))^l
\mbox{ for all } K \in \mathcal{T}_h \},
\end{align*}
where, $k=1,\;2$, $P^k$ denotes the polynomial approximation space, and $l=1,\;2$ indicates either a scalar or a vector.

In this paper, Generalized Multiscale Finite Element Method (GMsFEM) framework is applied.
In the GMsFEM methodology, one divides the computations into offline and online computations.
The offline computations are based on a preliminary dimension reduction of the
fine-grid finite element spaces (that may include
dealing with additionally important physical parameters, uncertainties and nonlinearities), and then the online procedure (if needed)
is applied to construct a reduced order model.
We start by constructing offline spaces.

\begin{figure}[tb]
  \centering
  \includegraphics[width=1.0 \textwidth]{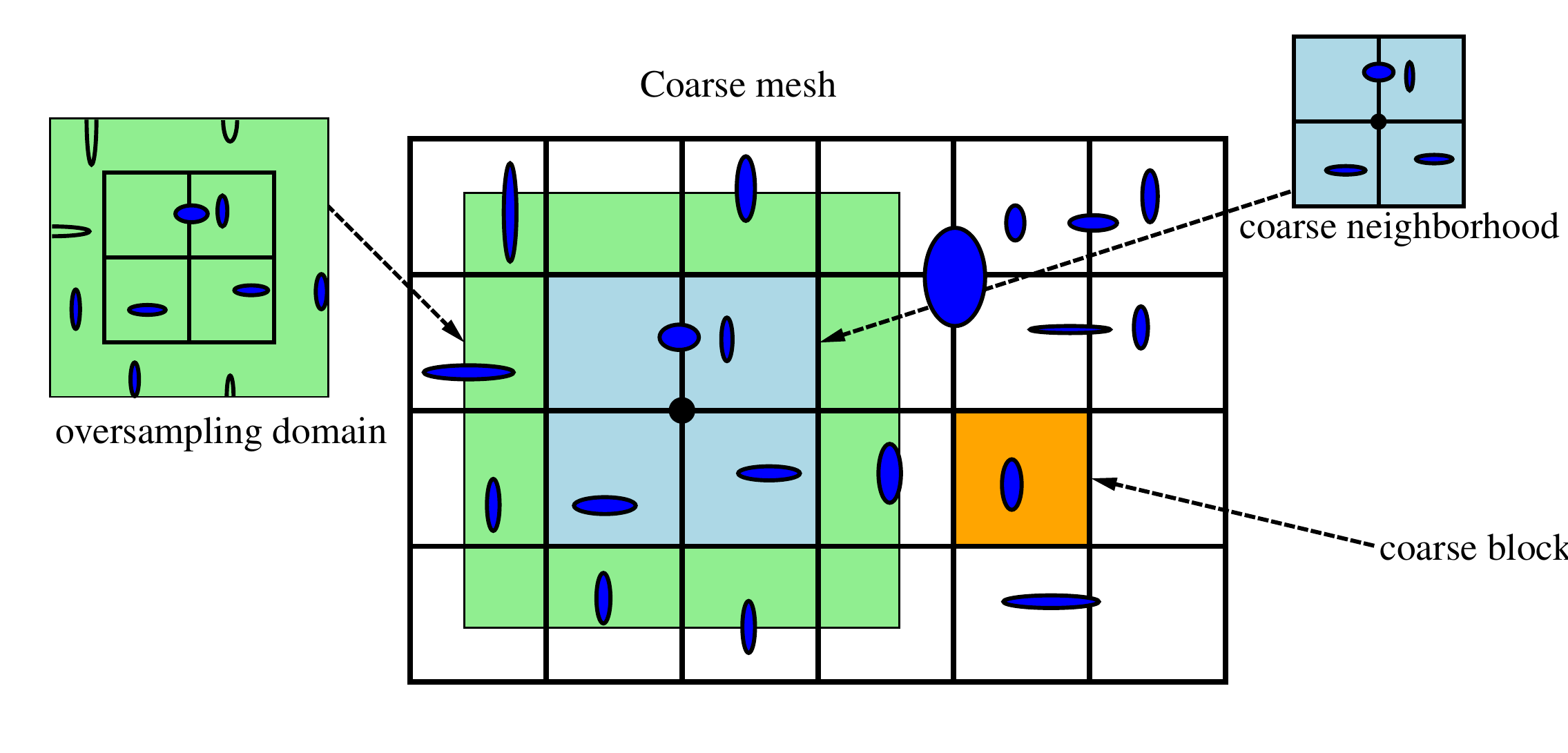}
  \caption{Illustration of a coarse neighborhood and oversampled domain}
  \label{schematic}
\end{figure}

\begin{figure}[tb]
  \centering
  \includegraphics[width=0.9 \textwidth]{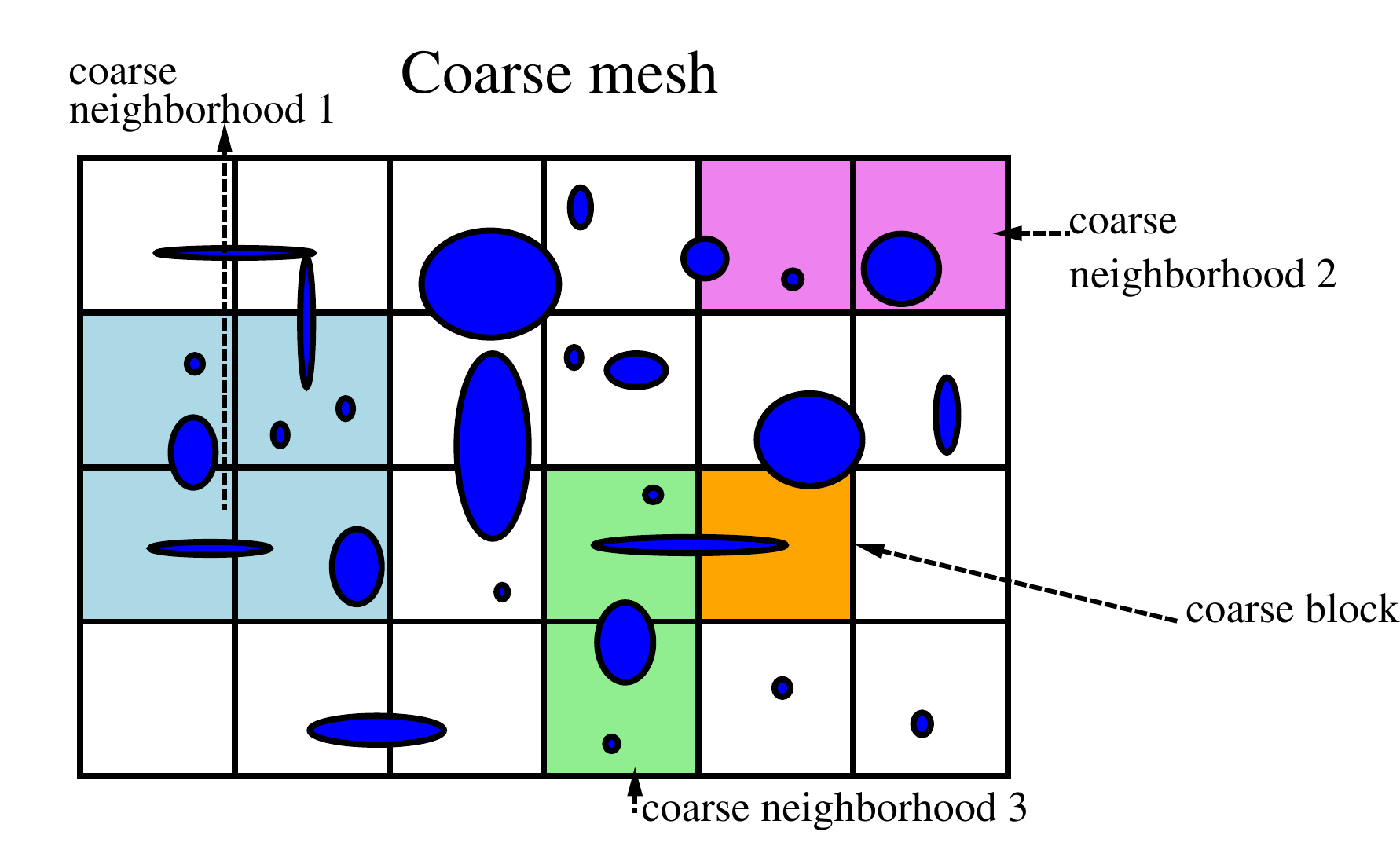}
  \caption{Illustration of a coarse neighborhood for Stokes equations}
  \label{schematic_stokes}
\end{figure}

We construct the coarse function space
$${V}_{\text{off}} := \mbox{span}\{{\phi}_i\}_{i=1}^{N},$$
where $N$ is the number of coarse basis functions.
Each ${\phi}_i$ is supported in some coarse
neighborhood $w_l$ (see Figures \ref{schematic} and \ref{schematic_stokes} for an illustration of coarse neighborhoods).

%The idea is then to work on the reduced spaces $V^\text{off}$ instead of the
%original spaces $V_h$.
%In the general GMsFEM methodology, these offline spaces are used  in the online
%computations where a further reduction may be  performed; see  \cite{egw10,egh12} for details.
The overall performance of the resulting GMsFEM depends on the
approximation properties of the resulting offline and online coarse spaces directly. Usually,
a spectral problem is involved for a good approximation of the local solution space.
In this paper, we focus on the construction of the offline spaces
since the differential operators $\mathcal {L}^{\eps}$
is parameter independent and linear.

The GMsFEM seeks an approximation ${u}_0\in{V}^{\text{off}}$,
 which satisfies the coarse-scale offline formulation,
\begin{align}\label{eq:coarse_system}
\avrg{\mathcal {L}^{\eps}(u_0),v}{\om^{\eps}}= (f,v)_{\om^{\eps}} \qquad  \text{for all } v \in V^{\text{off}}.
\end{align}
Recall that the definitions of the bilinear forms $\avrg{\mathcal {L}^{\eps}(u_0),v}{\om^{\eps}}$ are defined above, and $(f,v)_{\om^{\eps}}$ is the $L_2$ inner product.

We can interpret the method in the following way using matrix representations. Recall that the coarse basis functions
$\{{\phi}_i\}^{N_c}_{i=1}$ are defined on the fine grid, and
can be represented by the fine-grid basis functions. Specifically, we introduce the following matrices:
\[
R^T_0 = [{\phi}_1, \dots, {\phi}_{N_c}],
\]
where we identify the basis $\phi_i$ with their coefficient vectors
on the fine-grid basis.
Then, the matrix analogue of the system \eqref{eq:coarse_system} can be equivalently written as
\begin{align}\label{eqn:offlineSystem}
  R_0\avrg{\mathcal {L}^{\eps}(\phi_m),\phi_n}{\om^{\eps}}R_0^{T}u_0= R_0F,
\end{align}
where $F$ is the fine-grid discretization of $f$.
Further, once we solve the coarse system \eqref{eqn:offlineSystem}, we can recover the fine scale solution by $R_0^T u_0$.
In other words, $R^{T}_{0}$ can be regarded as the transformation (also known as interpolation or
downscaling) matrix  from the space ${V}_{\text{off}}$ to the space ${V}_{h}$.

The accuracy of the GMsFEM relies on the coarse basis functions $\{\phi_i\}$.
We shall present
the construction of suitable basis functions for the differential operators in Section \ref{locbasis}.
%------------------------------------
\section{Local basis functions}
\label{locbasis}
%------------------------------------
In this section we describe the offline-online computational procedure, and elaborate on some applicable choices for the associated bilinear forms to be used in the coarse space construction. Below, we offer a general outline for the procedure.

\begin{itemize}
\item[1.]  Offline computations:
\begin{itemize}
\item 1.0. Coarse grid generation.
\item 1.1. Construction of a snapshot space that will be used to compute an offline space.
\item 1.2. Construction of a small dimensional offline space by performing dimension reduction in the space of local snapshots.
\end{itemize}
\item[2.] Online computations:
\begin{itemize}
\item 2.1. For each input parameter, compute multiscale basis functions (it is needed for parameter-dependent or nonlinear problems).
\item 2.2. Solve a coarse-grid problem for a forcing term and boundary condition.
\item 2.3. Iterative solvers, if needed.
\end{itemize}
\end{itemize}

In the offline computation, we first construct a snapshot space $V_{\text{snap}}^{\omega_i^+}$ or $V_{\text{snap}}^{\omega_i}$, depending on the choice of domain to generate the snapshot space, where $\omega_i^+$ is an oversampled region that contains a coarse neighborhood $\omega_i$. Construction of the snapshot space involves solving the local problems for various choices of input parameters, and we describe the details below.

\subsection{Snapshot space}
\label{sec:harmonic}

First, we introduce the concept of coarse neighborhood as before. An illustration of coarse neighborhood is shown in Figures \ref{schematic} and \ref{schematic_stokes}, respectively, for Darcy  and Stokes equations.
The snapshot space is composed of harmonic extension
of fine-grid functions defined on the boundary of $\omega_i$ excluding the inactive cells in it. The local snapshots satisfy boundary conditions imposed on $\partial \mathcal {B}^{\eps}$.
More precisely, for each fine-grid function, $\delta_l^h(x)$,
which is defined by
$\delta_l^h(x)=\delta_{l,k},\,\forall l,k\in \textsl{J}_{h}(\omega_i)$, where $\textsl{J}_{h}(\omega_i)$ denotes the fine-grid boundary node on $\partial\omega_i\backslash \mathcal {B}^{\eps}$.

For parameter-independent problem, we
solve
\begin{align}\label{eq:harmonic}
\mathcal {L}^{\eps}(\psi_{l}^{ \text{snap}})=0\ \ \text{in} \ \omega_i\backslash \mathcal {B}^{\eps}
\end{align}
subject to the boundary condition
 $ \psi_{l}^{\text{snap}}=\delta_l^h(x)$ on $\partial\omega_i$ and
$ \psi_{l}^{\text{snap}}=0$ on $\partial\mathcal {B}^{\eps}\cap \bar{\omega_i}$.
Note that for the differential operator with Neumann boundary condition $ \displaystyle\frac{\partial u}{\partial n}=0$, we will use Neumann boundary condition instead for the local problems proposed above.

We denote the space composed of $\psi_{l}^{\text{snap}}$ as
$$
V_{\text{snap}}= \text{span}\{ \psi_{l}^{\text{snap}}: ~~ 1\leq l \leq M_i \},
$$
for each coarse neighborhood $\omega_i$, where $M_i$ denotes the number of snapshots in
the region $\omega_i$.
We emphasize that an oversampling strategy (refer to Figure \ref{schematic} for an illustration of oversampling domain)
can be applied for the construction of snapshots to obtain a higher accuracy (\cite{Efendiev_oversampling13}).
\begin{remark}
The snapshots described above are referred to as the harmonic basis in the numerical part.
\end{remark}

\begin{remark}
One can also use all the fine grid functions as snapshots.
These snapshots are referred to as the spectral basis in the numerical part.
\end{remark}

\subsection{Offline space}

This section is devoted to the construction of the offline space via a
 spectral decompostion.
In order to construct the  offline
space $V_{\text{off}}$, we
perform a dimension reduction in the snapshot space using an auxiliary
spectral decomposition. The main objective is to seek
a subspace of the snapshot space such that it can approximate any element of the snapshot space in the appropriate sense defined via auxiliary bilinear forms.
We will consider the following eigenvalue problems in the space of snapshots:
\begin{eqnarray}
A^{\text{off}} \Psi_k^{\text{off}} &=& \lambda_k^{\text{off}} S^{\text{off}}\Psi_k^{\text{off}}.  \label{offeig1}
%text{off}} &=& \lambda_k^{\text{off}} S^{+, \text{off}} \Psi_k^{\text{off}}
\end{eqnarray}
%  %Here, $A$ is constructed by integrating only on $\omega_i\backslash \mathcal {B}^{\eps}$: $A=\avrg{\mathcal {L}^{\eps}(\phi_m),\phi_n}{\omega_i\backslash \mathcal {B}^{\eps}}$, $\phi_m$ and $\phi_n$ are the $m^{th}$ and $n^{th}$ fine-scale basis over the domain $\omega_i\backslash \mathcal {B}^{\eps}$.
The definitions of $A^{\text{off}}$ and $S^{\text{off}}$ for Laplace, elasticity
 and Stokes equations are listed below.
  To generate the offline space we then choose the smallest $M_{\text{off}}$ eigenvalues from Eqn.~\eqref{offeig1} and form the corresponding eigenvectors in the respective
snapshot space  by setting
$\psi_k^{\text{off}} = \sum_j \Psi_{kj}^{\text{off}} \psi_j^{\text{snap}}$, for $k=1,\ldots, M_{\text{off}}$, where $\Psi_{kj}^{\text{off}}$ are the coordinates of the vector $\Psi_{k}^{\text{off}}$. We then create the offline matrices
 $$
R_{\text{off}} = \left[ \psi_{1}^{\text{off}}, \ldots, \psi_{M_{\text{off}}}^{\text{off}} \right].
$$

%\begin{example}\label{ex:eigen_prob}
\begin{itemize}
\item For Laplace operator \eqref{eqn:laplace},
\begin{equation}
\begin{split}
 A^{\text{off}}= [a^{\text{off}}_{mn}] = \int_{\omega_i \backslash\mathcal {B}^{\eps}}\nabla\psi_m^{\text{snap}}\cdot\nabla\psi_n^{\text{snap}}, \text{ and }\\
 S^{\text{off}}= [s^{\text{off}}_{mn}] =\int_{\omega_i \backslash\mathcal {B}^{\eps}}\psi_m^{\text{snap}}\psi_n^{\text{snap}}.
\end{split}
\end{equation}
\item For elasticity operator \eqref{eqn:elasticity},
\begin{equation}
\begin{split}
 A^{\text{off}}= [a^{\text{off}}_{mn}] =
\int_{\omega_i \backslash\mathcal {B}^{\eps}}
\Big( 2\mu  {\varepsilon}( \psi_m^{\text{snap}}) :  {\varepsilon}( \psi_n^{\text{snap}})
+ \lambda \nabla\cdot  {\psi_m^{\text{snap}}} \, \nabla\cdot  \psi_n^{\text{snap}}
\Big) \;, \text{ and }\\
\displaystyle S^{\text{off}}= [s^{\text{off}}_{mn}] = \int_{\omega_i \backslash\mathcal {B}^{\eps}}  (\lambda + 2 \mu)  \psi_m^{\text{snap}} \cdot  \psi_n^{\text{snap}}.
\end{split}
\end{equation}

\item For Stokes operator \eqref{eqn:stokes},
\begin{equation}
\begin{split}
 A^{\text{off}}= [a^{\text{off}}_{mn}] = \int_{\omega_i \backslash\mathcal {B}^{\eps}}\nabla\psi_m^{\text{snap}}:\nabla\psi_n^{\text{snap}}, \text{\text{and}}\\
\displaystyle S^{\text{off}}= [s^{\text{off}}_{mn}] =\int_{\omega_i \backslash\mathcal {B}^{\eps}}\psi_m^{\text{snap}}\cdot\psi_n^{\text{snap}}.
\end{split}
\end{equation}
Note that for the Stokes operator, multiscale spaces are constructed for velocity field only, and
we use piecewise constant approximation on the coarse grid for the pressure.
\end{itemize}
%\end{example}
%
%\begin{remark}
%The spectral problems proposed above can be used to generate snapshot spaces as well (\cite{GMsFEM13, Efendiev_oversampling13}), which are referred to as spectral basis later.
%\end{remark}
%to be used in the online space construction.

We remark that some adaptive procedures for choosing the offline basis are proposed and analyzed in \cite{AdaptiveGMsFEM,AdaptiveGMsDGM}. In our future work,
we plan to present an analysis of the method which can guide us in choosing
the local spectral problems.

\section{Simulation results}
\label{sec:numerical}

\begin{figure}
\begin{center}
\begin{tabular}{c c}
\includegraphics[width=0.4\linewidth]{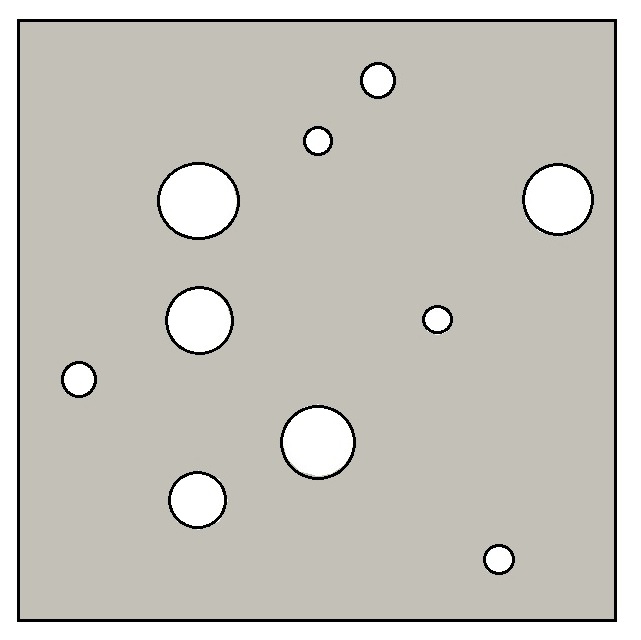}&\includegraphics[width=0.4\linewidth]{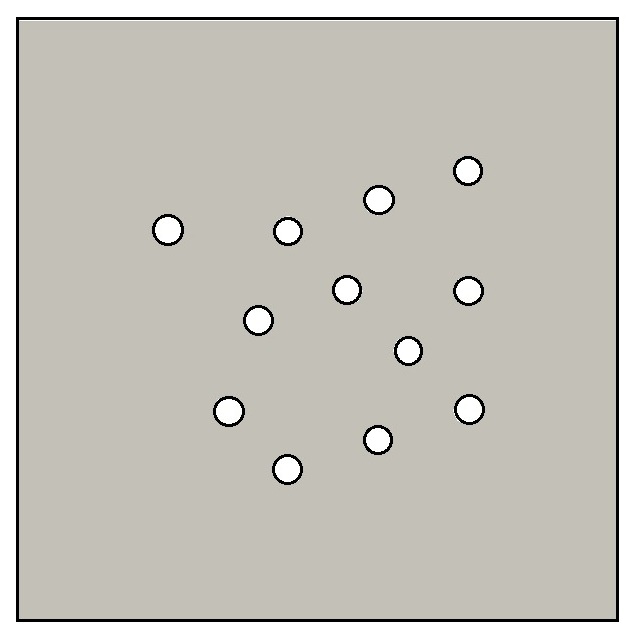}\\
%{Computational domain 1 (Big holes)}&{Computational domain 2 (Small holes)}
\label{fig:domain}
\end{tabular}
\caption{Two heterogeneous perforated medium used in the simulations.}
\end{center}
\end{figure}

\begin{figure}
\begin{center}
\includegraphics[width=0.8\linewidth]{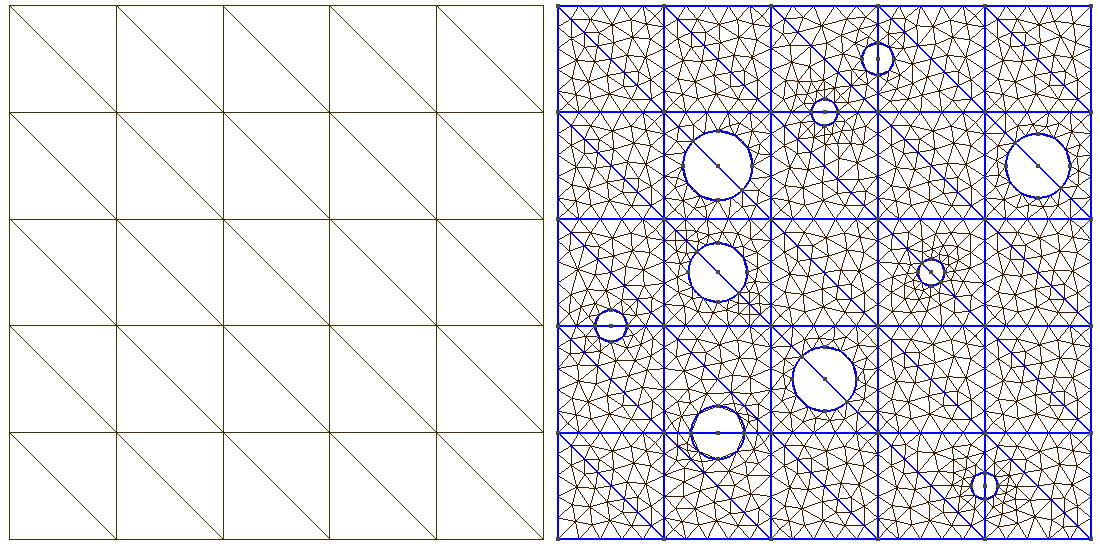}
\caption{Coarse-scale and fine-scale discretization corresponds to the heterogeneous medium on the left of Figure \ref{fig:domain}.}
\label{fig:mesh1}
\end{center}
\end{figure}

\begin{figure}
\begin{center}
\includegraphics[width=0.8\linewidth]{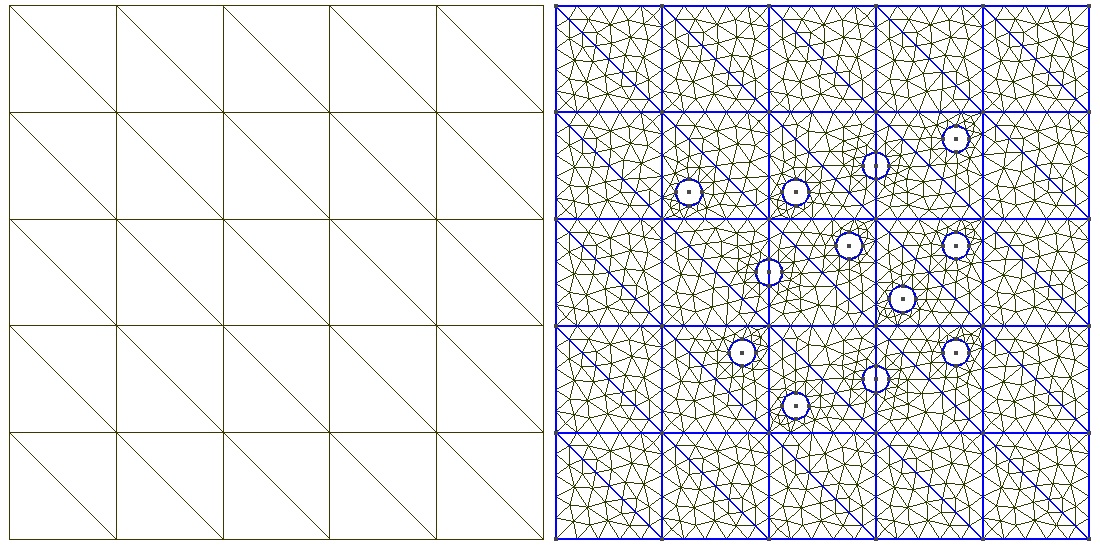}
\caption{Coarse-scale and fine-scale discretization corresponds to the  heterogeneous medium on the right of Figure \ref{fig:domain}.}
\label{fig:mesh2}
\end{center}
\end{figure}

In this section, we present simulation
results using the framework presented in Section \ref{locbasis}
for Laplace equation, elasticity equation and Stokes equations,
respectively. We set $\Omega = [0,1] \times [0,1]$ and use two types of perforated domains as illustrated in Figure \ref{fig:domain}, where the perforated regions $\mathcal {B}^{\eps}$ are circular. Note that we can also use perforated regions of other shapes instead and obtained similar results.  The computational domain is discretized coarsely using uniform triangulation as shown in the left of Figures \ref{fig:mesh1} and \ref{fig:mesh2}, where the coarse mesh size $H={1\over 5}$. Furthermore, nonuniform triangulation is used inside each coarse triangle element to obtain a finer discretization.
Examples of this triangulation are
displayed on the right of Figures \ref{fig:mesh1} and \ref{fig:mesh2}, where $2188$ and $2364$ nonuniform triangle elements are generated, respectively.


\subsection{Laplace equation in perforated domain}

First, we consider the Laplace operator \eqref{eqn:laplace}
imposed with zero Dirichlet boundary condition on the holes
$\partial \Omega^{\eps}\cap  \partial \mathcal {B}^{\eps}$
and $u=1$ on $\partial \Omega$, and $f=0$.

\begin{figure}[htp]
\begin{center}
\includegraphics[width=0.78\linewidth]{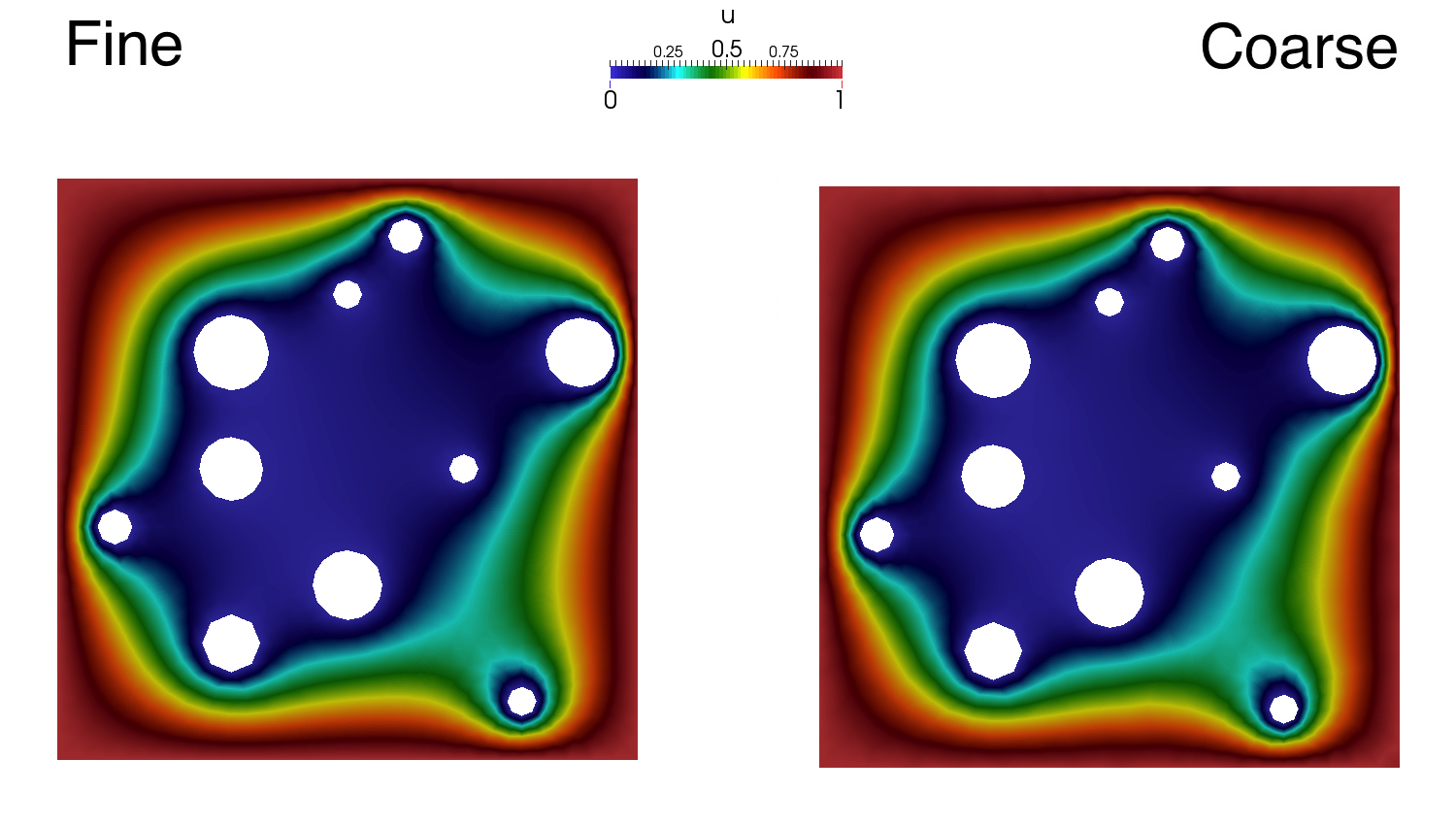}
\caption{Fine-scale and coarse-scale solution to the Laplace equation in the heterogeneous medium on the left of Figure \ref{fig:domain}. The dimension of the coarse space is 432.}
\label{fig:la1}
\end{center}
\end{figure}

\begin{figure}[htp]
\begin{center}
\includegraphics[width=0.78\linewidth]{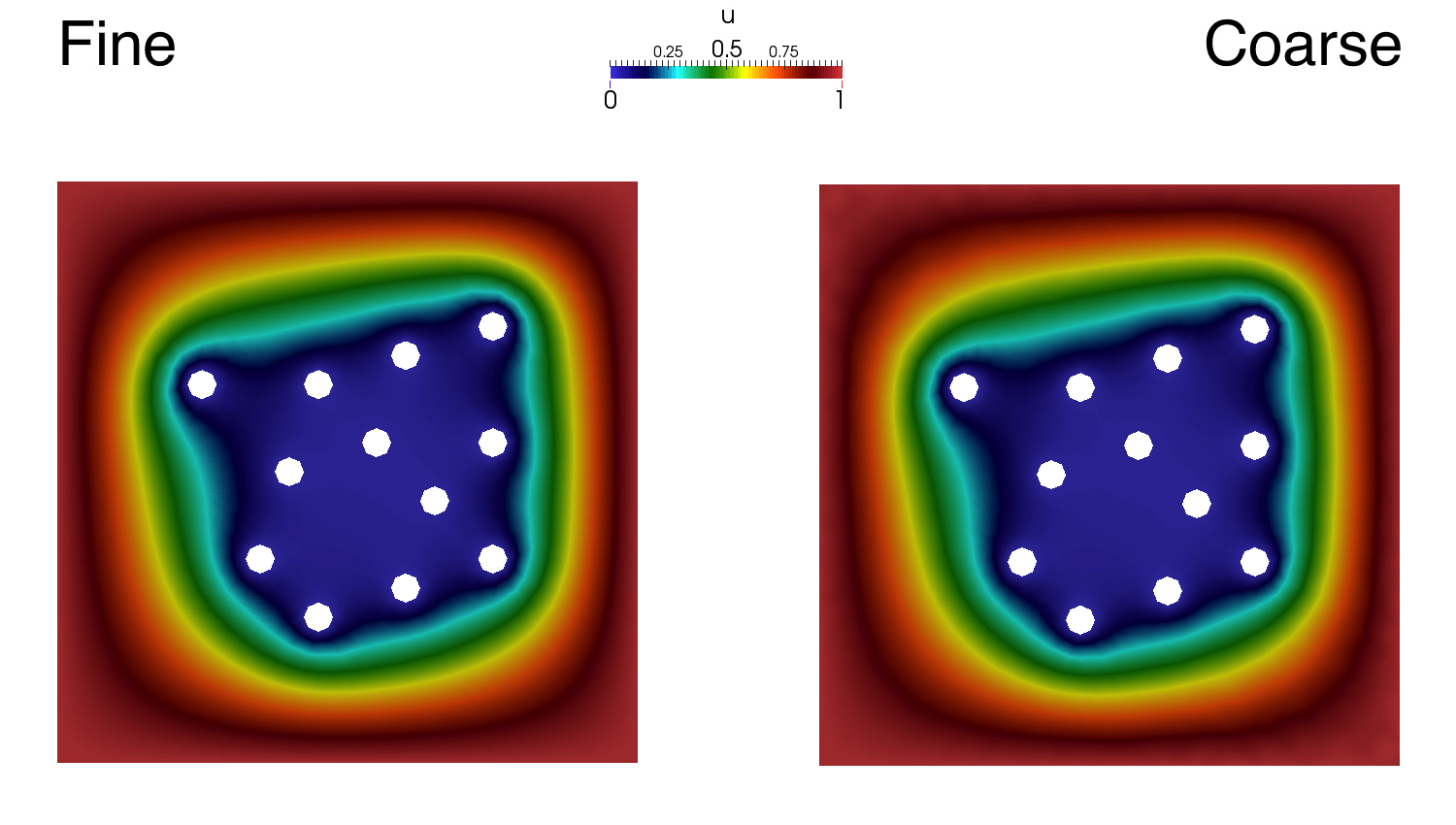}
\caption{Fine-scale and coarse-scale solution to the Laplace equation in the heterogeneous medium on the right of Figure \ref{fig:domain}. The dimension of the coarse space is 432.}
\label{fig:la2}
\end{center}
\end{figure}

The simulation results in perforated
domains as shown in Figure \ref{fig:domain} are
illustrated in Figures \ref{fig:la1} and \ref{fig:la2},
respectively.
The multiscale solution is obtained in an offline space of dimension $432$ (using $12$ basis functions per coarse neighborhood)
and the fine-scale reference solution is obtained in a space of dimension $1187$.
Compared the fine-scale solution on the left with
 the coarse-scale solution on the right of the figures, we observe that
the GMsFEM can approximate the fine-scale solution accurately.


Besides, two types of snapshot spaces are applied in our
simulations with the results shown in Tables \ref{tab:la1} and \ref{tab:la2}.
The first column of each table shows
 the number of basis in each coarse node ($N_c$).
The dimensions of the offline spaces are given in the second column.
 The next two columns display the $L_2$ and $H_1$ relative errors
for  using spectral basis.
We can see that the $L_2$ and $H_1$ relative errors are $0.2\%$ and $2\%$ respectively when the dimension of the coarse space is $432$.
The results applying harmonic extenstion type of
basis (Eqn. \eqref{eq:harmonic}) are listed in the last two
columns.
As displayed in Table \ref{tab:la1}, the $L_2$ and $H_1$ relative errors are $1\%$ and $17\%$ respectively when the dimension of the coarse space is $432$.
%Here, the dimension of the fine-scale solution is 1187.
Note that the GMsFEM gives a very good solution in both cases with only about $36\%$
of unknown compared to the fine-scale solver.
We also observe in Tables \ref{tab:la1} and \ref{tab:la2}
that a  rapid decay in the errors as we increase the number of basis functions.

\begin{table}
\centering
\begin{tabular}[hp]{|c|c|cc|cc|}
\hline
\multirow{2}{*}{$N_c$} & \multirow{2}{*}{dim} & \multicolumn{2}{c|}{harmonic basis} & \multicolumn{2}{c|}{spectral basis} \\
 &  & $L_2$ & $H_1$ & $L_2$ & $H_1$ \\
\hline \hline
1   & 36   		& 0.18   & 0.75  & 0.19 & 0.63\\
2   & 72 		& 0.12   & 0.65  & 0.10 & 0.45\\
4   & 144 		& 0.06   & 0.48  & 0.02 & 0.18\\
8   & 288 		& 0.03   & 0.35  & 0.01 & 0.09\\
12  & 432 		& 0.01   & 0.17  & 0.002 & 0.02\\
16  & 576  	    & 0.004  & 0.04  & 0.001 & 0.009\\
\hline
\end{tabular}
\caption{Numerical tests for Laplace operator in heterogeneous medium shown on the left of Figure \ref{fig:domain}. Fine-scale problem dimension is 1187.}
\label{tab:la1}
\end{table}
\begin{table}
\begin{center}
\begin{tabular}[hp]{|c|c|cc|cc|}
\hline
\multirow{2}{*}{$N_c$} & \multirow{2}{*}{dim} & \multicolumn{2}{c|}{harmonic basis} & \multicolumn{2}{c|}{spectral basis} \\
 &  & $L_2$ & $H_1$ & $L_2$ & $H_1$ \\
\hline \hline
1   & 36   		& 0.10   & 0.89  & 0.06     & 0.34\\
2   & 72 		& 0.09   & 0.88  & 0.03     & 0.25\\
4   & 144 		& 0.05   & 0.65  & 0.01     & 0.10\\
8   & 288 		& 0.03   & 0.52  & 0.005   & 0.05\\
12 & 432 		& 0.01   & 0.29  & 0.001   & 0.02\\
16 & 576  	& 0.003 & 0.06  & 0.0006 & 0.009\\
\hline
\end{tabular}
\end{center}
\caption{Numerical tests for Laplace operator in heterogeneous medium shown on the right of Figure \ref{fig:domain}. Fine-scale problem dimension is 1269.}
\label{tab:la2}
\end{table}

\subsection{Elasticity equation in perforated domain}

\begin{figure}
\begin{center}
\includegraphics[width=0.65\linewidth]{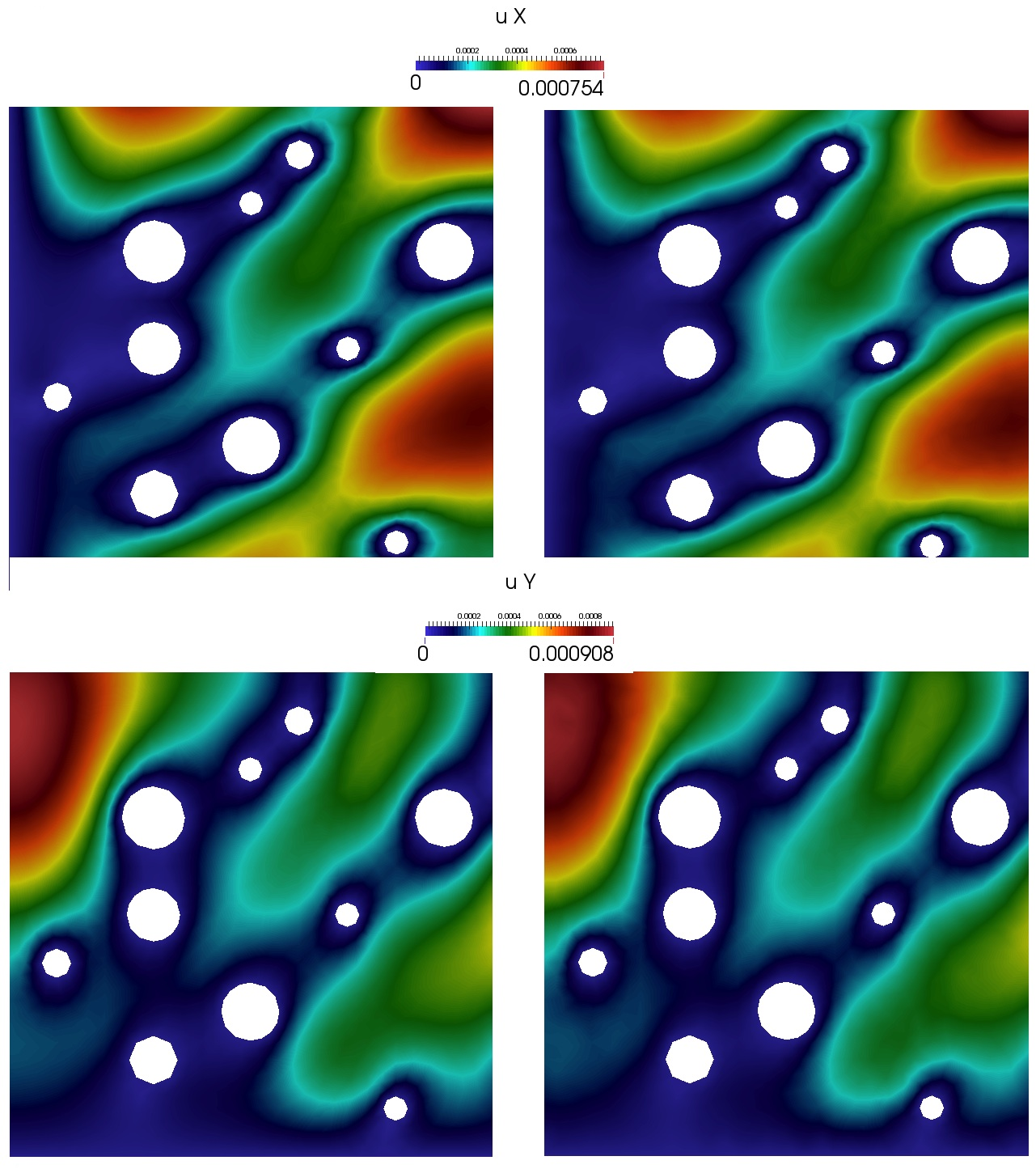}
\caption{The fine-scale and coarse-scale solutions of the elasticity equation correspond to the heterogeneous perforated domain shown on the left of Figure \ref{fig:domain}. The dimension of the fine-scale solution is 2374 and the dimension of the coarse space is 864.}
\label{fig:el1}
\end{center}
\end{figure}

\begin{figure}
\begin{center}
\includegraphics[width=0.65\linewidth]{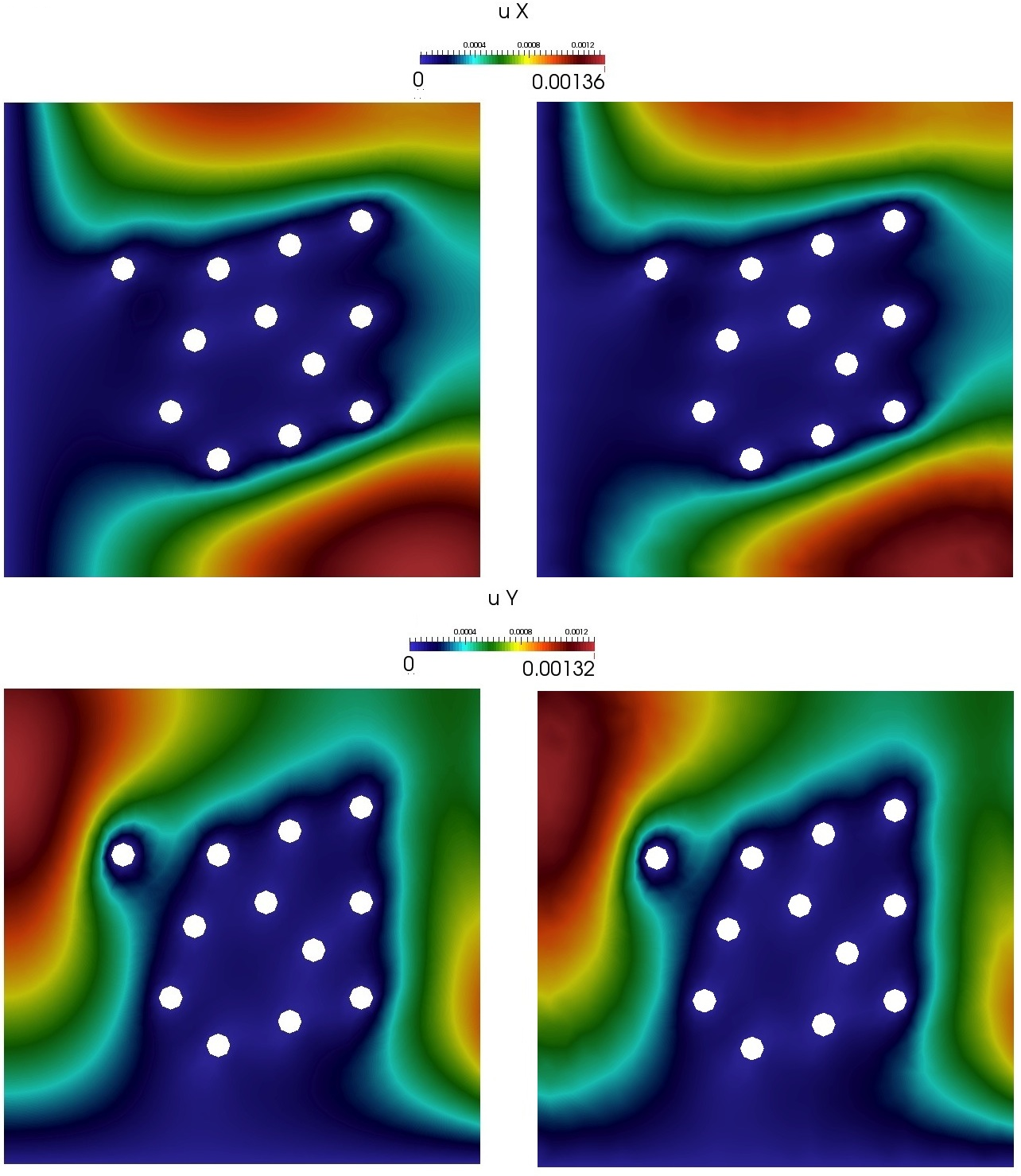}
\caption{The fine-scale and coarse-scale solutions of the elasticity equation correspond to the heterogeneous perforated domain shown on the right of Figure \ref{fig:domain}. The dimension of the fine-scale solution is 2538 and the dimension of the coarse space is 900.}
\label{fig:el2}
\end{center}
\end{figure}

Next, we consider elasticity operator \eqref{eqn:elasticity}.
We use zero displacements
$u =0$ on the holes,
$u_x = 0, \sigma_y = 0$ on the left boundary,
$\sigma_x = 0, u_y = 0$  on the bottom boundary and
$\sigma_x = 0, \sigma_y = 0$ on the top and right boundaries.
Here, $u = (u_x, u_y) $ and $\sigma = (\sigma_x, \sigma_y) $.
The source term is defined by $f = (10^7, 10^7)$,
the elastic modulus is given by $E = 10^9$, Poisson's ratio is $\nu = 0.22$, where
\[
\mu = \frac{E}{2 (1 + \nu)}, \quad
\lambda = \frac{E \nu}{(1+ \nu) ( 1- 2 \nu)}.
\]
The fine-scale solution and coarse-scale solution corresponding to different perforated domains in Figure \ref{fig:domain} are presented in Figures \ref{fig:el1} and \ref{fig:el2}. The fine-scale displacement is displayed on the left, while the coarse-scale displacement is displayed on the right.
The multiscale solution is obtained in an offline space of dimension $864$ (using $12$ basis functions per coarse neighborhood for each component of the displacement)
and the fine-scale reference solution is obtained in a space of dimension $2374$.
Comparing the fine-scale solution with the coarse-scale solution in Figures \ref{fig:el1} and \ref{fig:el2}, we can observe a good accuracy.

Furthermore, the relative errors for different dimensions of coarse spaces with harmonic basis functions (Eqn. \eqref{eq:harmonic}) and spectral basis functions are shown in Tables \ref{tab:el1} and \ref{tab:el2}, respectively. As displayed in Table \ref{tab:el1}, the $L_2$ and $H_1$ relative errors using harmonic basis
are $3\%$ and $14\%$, respectively, when the dimension of the coarse space is $1152$ for the perforated domain shown in the left of Figure \ref{fig:domain}.
In addition, for the same dimension of the coarse space and the same domain, the $L_2$ and $H_1$ relative errors using spectral basis
are $0.5\%$ and $5\%$, respectively
Note that the dimension of the fine-scale solution is 2374.
Thus, we see that the proposed GMsFEM gives accurate solutions.
We also observe that GMsFEM gives more accurate solution
when more basis are included in the offline space.
We remark that the results (see Table \ref{tab:el2})
follow a similar pattern for the domain shown on the right
of Figure \ref{fig:domain}.

\begin{table}
\begin{center}
\begin{tabular}[hp]{|c|c|cc|cc|}
\hline
\multirow{2}{*}{$N_c$} & \multirow{2}{*}{dim} & \multicolumn{2}{c|}{harmonic basis} & \multicolumn{2}{c|}{spectral basis} \\
 &  & $L_2$ & $H_1$ & $L_2$ & $H_1$ \\
\hline \hline
1   & 72   		& 0.65 & 0.74  & 0.58 & 0.75\\
2   & 144 		& 0.41 & 0.57  & 0.21 & 0.46\\
4   & 288 		& 0.32 & 0.48  & 0.07 & 0.26\\
8   & 576 		& 0.16 & 0.32  & 0.02 & 0.13\\
12 & 864 		& 0.11 & 0.25  & 0.009 & 0.08\\
16 & 1152 	& 0.03 & 0.14 & 0.005 & 0.05\\
20 & 1440 	& 0.01 & 0.04 & 0.003 & 0.03\\
\hline
\end{tabular}
\end{center}
\caption{Numerical tests for elasticity equations in heterogeneous medium shown on the left of Figure \ref{fig:domain}. {The fine-scale problem dimension is 2374.}}
\label{tab:el1}
\end{table}

\begin{table}
\begin{center}
\begin{tabular}[hp]{|c|c|cc|cc|}
\hline
\multirow{2}{*}{$N_c$} & \multirow{2}{*}{dim} & \multicolumn{2}{c|}{harmonic basis} & \multicolumn{2}{c|}{spectral basis} \\
 &  & $L_2$ & $H_1$ & $L_2$ & $H_1$ \\
\hline \hline
1   & 72   		& 0.79 & 0.86  & 0.65 & 0.81  \\
2   & 144 		& 0.53 & 0.69  & 0.11 & 0.35  \\
4   & 288 		& 0.45 & 0.62  & 0.04 & 0.20  \\
8   & 576 		& 0.22 & 0.41  & 0.01 & 0.12  \\
12 & 864 		& 0.16 & 0.34  & 0.007 & 0.07  \\
16 & 1152 	& 0.06 & 0.20  & 0.003 & 0.05 \\
20 & 1440 	& 0.01 & 0.12  & 0.002 & 0.04 \\
\hline
\end{tabular}
\end{center}
\caption{Numerical tests for elasticity equations in heterogeneous medium shown on the right of Figure \ref{fig:domain}.{ The fine-scale problem dimension is 2538.} }
\label{tab:el2}
\end{table}

\subsection{Stokes equations in perforated domain}

\begin{figure}[htp]
\begin{center}
\includegraphics[width=0.65\linewidth]{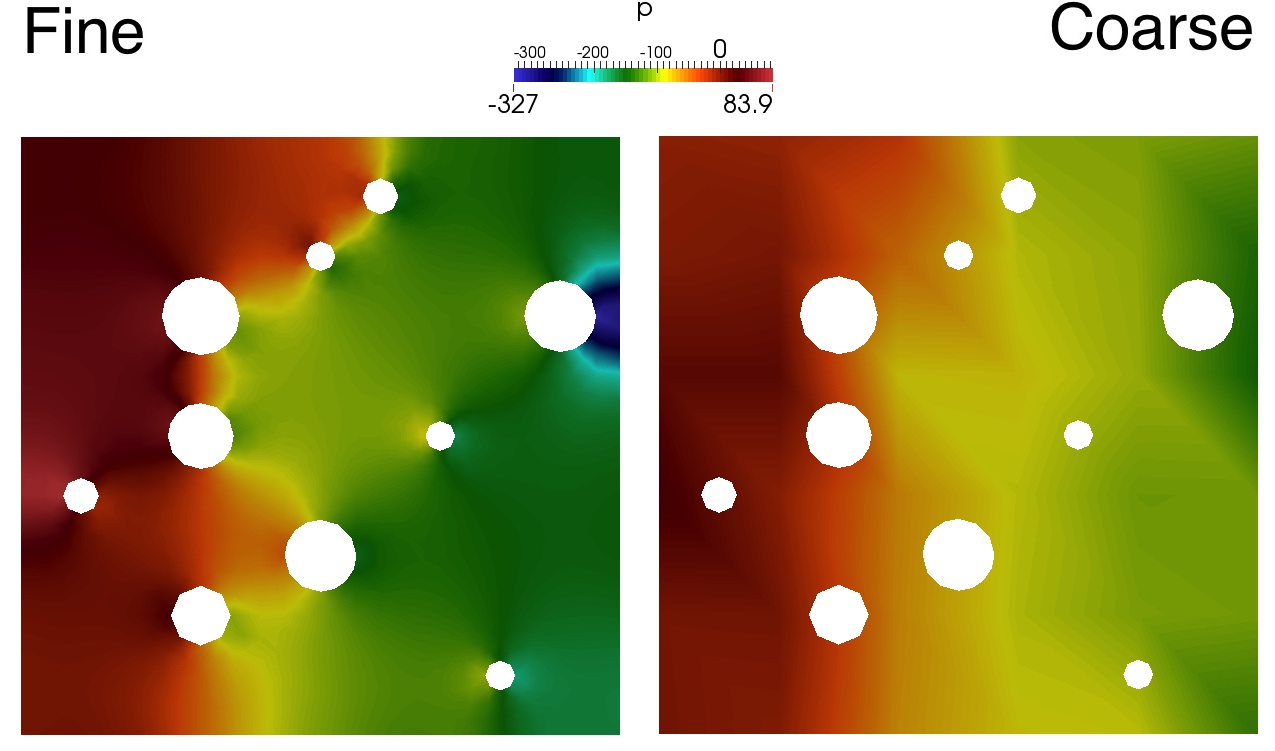}
\includegraphics[width=0.65\linewidth]{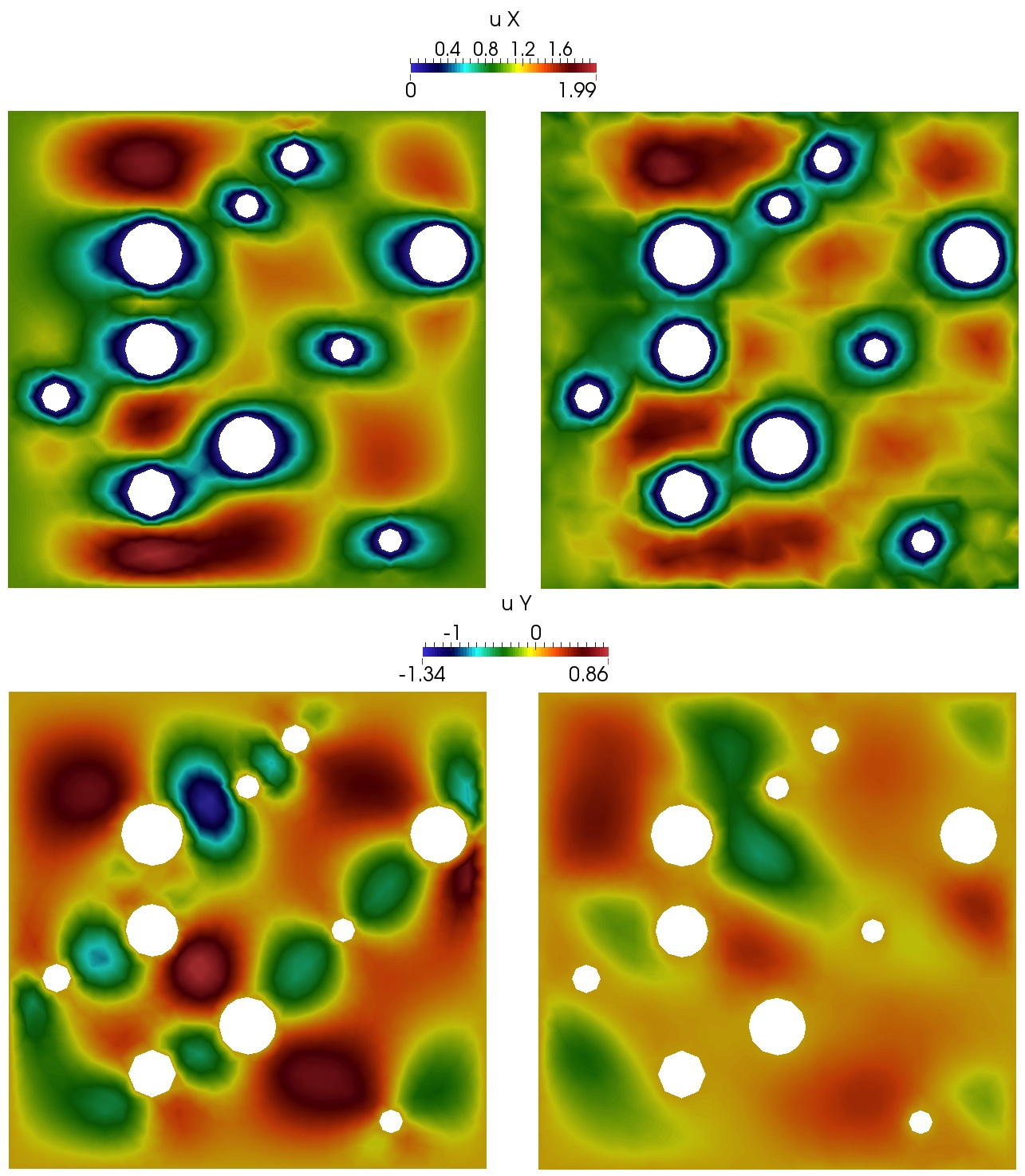}
\caption{The fine-scale and coarse-scale solutions of the pressure, x-component and y-component of the velocity correspond to the heterogeneous perforated domain shown on the left of Figure \ref{fig:domain}. The dimension of the fine-scale solution is 3561.}
\label{fig:st11}
\end{center}
\end{figure}

\begin{figure}[htp]
\begin{center}
\includegraphics[width=0.65\linewidth]{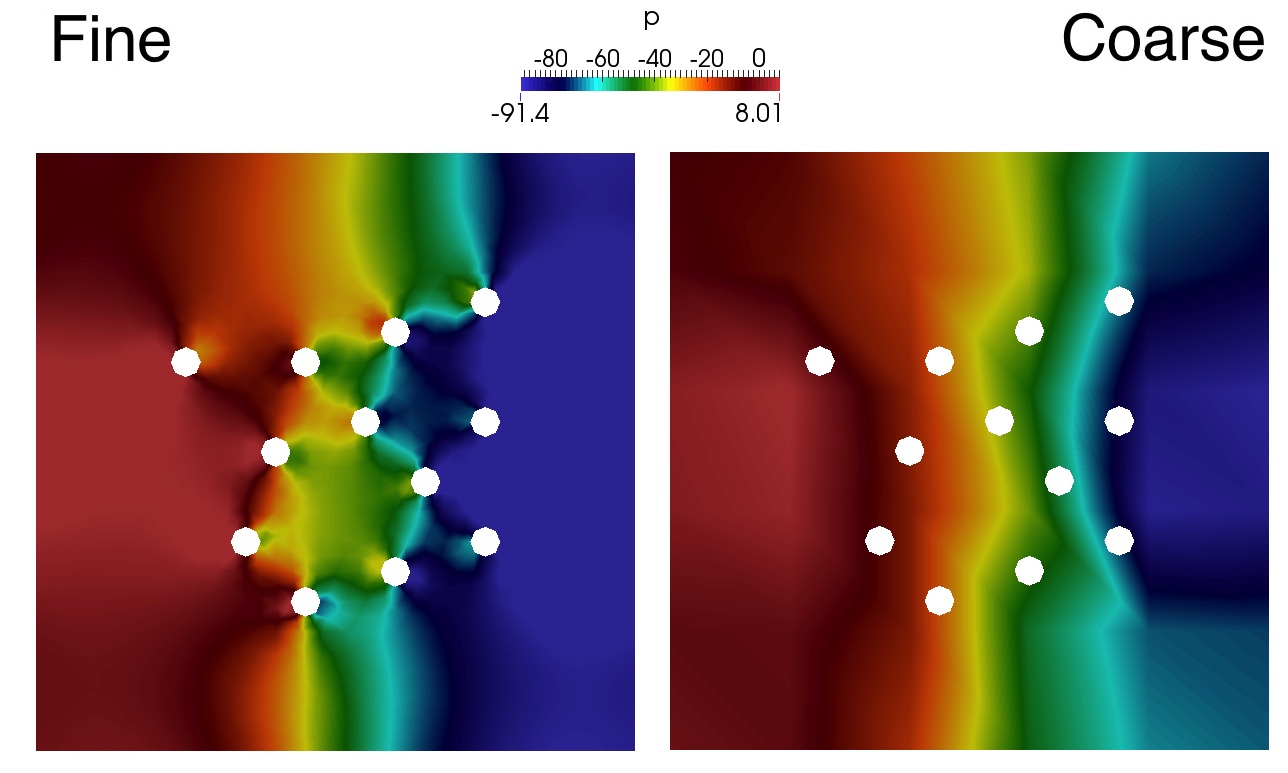}
\includegraphics[width=0.65\linewidth]{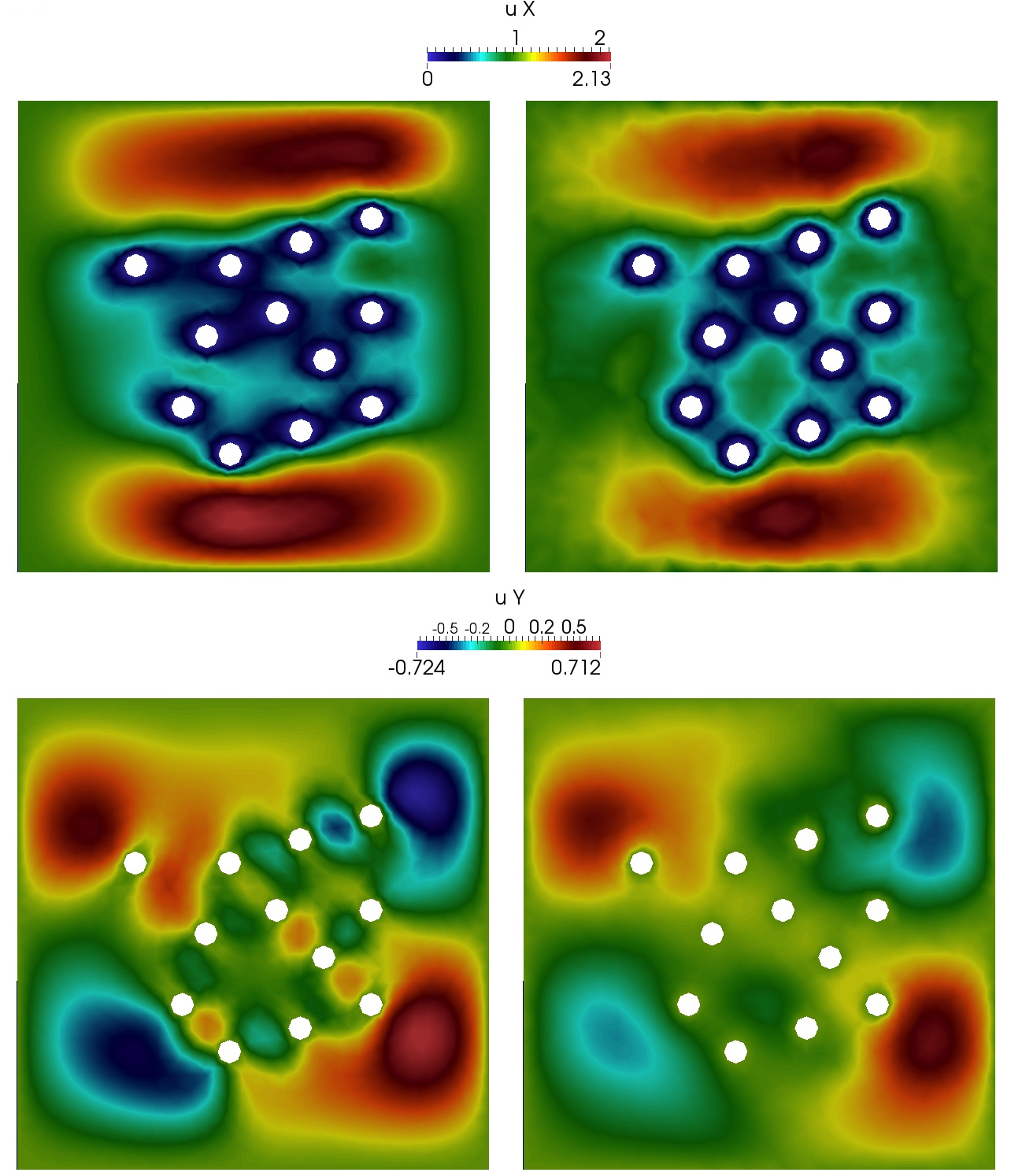}
\caption{The fine-scale and coarse-scale solutions of the pressure, x-component and y-component of the velocity correspond to the heterogeneous perforated domain shown on the right of Figure \ref{fig:domain}. The dimension of the fine-scale solution is 3807.}
\label{fig:st21}
\end{center}
\end{figure}

In our final example, we consider the Stokes operator \eqref{eqn:stokes}
with zero velocity $u = (0, 0)$ on $\partial \Omega^{\eps}\cap  \partial \mathcal {B}^{\eps}$ and $u = (1, 0)$ on $\partial \Omega$.
The numerical results corresponding to different perforated domains in Figure \ref{fig:domain} are listed in Figures \ref{fig:st11} and \ref{fig:st21}. The fine-scale solution and coarse-scale solution are depicted on the left and right of Figures \ref{fig:st11} and \ref{fig:st21}, respectively. To improve the accuracy,
we have enriched pressure spaces by considering a splitting algorithm
\cite{chorin1968numerical}. Our preliminary numerical results show
an improvement.

\begin{table}[htp]
\begin{center}
\begin{tabular}[hp]{|c|c|cc|}
\hline
\multirow{2}{*}{$N_c$} & \multirow{2}{*}{dim} & \multicolumn{2}{c|}{harmonic basis} \\% & \multicolumn{2}{c|}{spectral problem 2} \\
 &  & $L_2$ & $H_1$ \\% & $L_2$ & $H$ \\
\hline \hline
1   & 108  & 0.85 & 0.88 \\%  & 0.86 & 0.97\\
2   & 180 	& 0.68 & 0.78  \\% & 0.49 & 0.76\\
4   & 324 	& 0.40 & 0.69  \\% & 0.33 & 0.65\\
8   & 612 	& 0.26 & 0.67  \\% & 0.19 & 0.56\\
12  & 900  	& 0.23 & 0.60  \\% & 0.17 & 0.50\\
16  & 1188 	& 0.23 & 0.58  \\% & 0.17 & 0.49\\
\hline
\end{tabular}
\end{center}
\caption{Numerical tests for Stokes operator in heterogeneous medium shown on the left of Figure \ref{fig:domain}. Fine-scale problem dimension is 3561.} %. Fine-scale dim 3561. For $\bm g = (1, 0)$ and $\bm g = (1, 1)$}
\label{tab:st1}
\end{table}
\begin{table}[htp]
\begin{center}
\begin{tabular}[hp]{|c|c|cc|}
\hline
\multirow{2}{*}{$N_c$} & \multirow{2}{*}{dim} & \multicolumn{2}{c|}{harmonic basis} \\%  & \multicolumn{2}{c|}{spectral problem 2} \\
 &  & $L_2$ & $H_1$\\%  & $L_2$ & $H$ \\
\hline \hline
1   & 108  	& 0.65 & 0.87  \\% & 0.79 & 0.92\\
2   & 180 		& 0.56 & 0.82  \\% & 0.53 & 0.86\\
4   & 324 		& 0.31 & 0.71  \\% & 0.30 & 0.71\\
8   & 612 		& 0.16 & 0.62  \\% & 0.14 & 0.61\\
12 & 900   	& 0.13 & 0.46  \\% & 0.11 & 0.41\\
16 & 1188 	& 0.13 & 0.44  \\% & 0.11 & 0.41\\
\hline
\end{tabular}
\end{center}
\caption{Numerical tests for Stokes operator in heterogeneous medium shown on the right of Figure \ref{fig:domain}. Fine-scale problem dimension is 3807.} %. Fine-scale dim 3807.  For $\bm g = (1, 0)$ and $\bm g = (1, 1)$}
\label{tab:st2}
\end{table}
The $L_2$ and $H_1$ relative errors in different perforated domains in Figure \ref{fig:domain} are listed in Tables
\ref{tab:st1} and \ref{tab:st2}, respectively. We observe that the $L_2$ and $H_1$ relative errors are $23\%$ and $58\%$, respectively. As more basis are taken in the construction of the offline space, the $L_2$ and $H_1$ relative errors become smaller.  In general, enriching pressure spaces and
constructing associated velocity spaces can improve the accuracy.
This is currently under investigation.

%\begin{figure}[htp]
%\begin{center}
%\includegraphics[width=0.7\linewidth]{stokes1/p11}
%\includegraphics[width=0.65\linewidth]{stokes1/u11}
%\caption{Pressure and velocity ($\bm g = (1, 1)$). Case 1 (Big holes)}
%\label{fig:st12}
%\end{center}
%\end{figure}

%\begin{figure}[htp]
%\begin{center}
%\includegraphics[width=0.7\linewidth]{stokes2/p11}
%\includegraphics[width=0.65\linewidth]{stokes2/u11}
%\caption{Pressure and velocity ($\bm g = (1, 1)$). Case 2 (Small holes)}
%\label{fig:st22}
%\end{center}
%\end{figure}

\subsection{Randomized snapshots for GMsFEM}

In this subsection, we will investigate the oversampling randomized
algorithm proposed in
\cite{cegl14} (shown in Table \ref{algorithm:random_snapshots11}).
The advantage of this algorithm lies in the fact that a much fewer number of snapshot basis functions are calculated maintaining a good accuracy to the solution space concurrently. Besides, the oversampling strategy is used to reduce the mismatching effects of boundary conditions imposed artificially in the construction of snapshot basis functions. For the sake of brevity, the results for elasticity equation are shown only.

The simulation results are presented in Tables \ref{tab:ov1} and \ref{tab:ov2} for elasticity problem. The results without oversampling are shown on the top of each table, while those with the oversampling strategy are shown on the bottom part of each table. In our simulation, we set the oversampling size $t=2$ (i.e., two extra fine-grid blocks are added to the original coarse region) and a buffer number $p_{\text{bf}}^{\omega_i}=4$ for each coarse neighborhood $\omega_i$, i.e., we will generate 4 more snapshot basis functions in each coarse neighborhood. In the numerical results, we report the fraction of the snapshots computed compared to all snapshot vectors. For the sake of completeness,
we list the algorithm in Table \ref{algorithm:random_snapshots11}.

\begin{table}[htp]
 \caption{Randomized GMsFEM Algorithm}
 \begin{tabular}{r l}
 \hline\hline
 \\
    \textbf{Input}:& Fine grid size $h$, coarse grid size $H$,
oversampling size $t$, buffer number $p_{\text{bf}}^{\omega_i}$ for each $\omega_i$, \\
    & the number of local basis functions $k_{\text{nb}}^{\omega_i}$ for each $\omega_i$;\\
    \textbf{output}: & Coarse-scale solution $u_H$.\\

  1.& Generate oversampling region for each coarse block: $\mathcal{T}^H$, $\mathcal{T}^h$, and $\omega_i^{+}$; \\

  2.& Generate $k_{\text{nb}}^{\omega_i}+p_{\text{bf}}^{\omega_i}$ random vectors $r_l$ and obtain randomized snapshots in $\omega_i^{+}$;  \\

& Add a snapshot that represents the constant function on $\omega_i^+$;\\
  3. & Obtain $k_{\text{nb}}^{\omega_i}$ offline basis by a spectral decomposition restricted to the original region; \\

4. & Construct multiscale basis functions and solve it.\\

 \hline\hline
 \end{tabular}
\label{algorithm:random_snapshots11}
 \end{table}

\begin{table}
\centering
\begin{tabular}[hp]{|c|cc|cc|}
\hline
\multirow{2}{*}{$N_c$}  &
\multicolumn{2}{c|}{full snapshots} &
\multicolumn{2}{c|}{randomized snapshots} \\
 & $L_2$ & $H_1$ & $L_2$ & $H_1$ \\
 \hline \hline
\multicolumn{5}{|c|}{without oversampling, $w_i$} \\ \hline
 &
\multicolumn{2}{c|}{100 \%} &
\multicolumn{2}{c|}{39.7 \%} \\  \hline
1 	 & 0.572 & 0.753 & 0.764 & 0.853 \\
2 	 & 0.217 & 0.466 & 0.592 & 0.742 \\
4 	 & 0.071 & 0.261 & 0.317 & 0.529 \\
8 	 & 0.023 & 0.136 & 0.173 & 0.386 \\
12 & 0.009 & 0.079 & 0.101 & 0.286 \\
16 & 0.005 & 0.054 & 0.055 & 0.214 \\
\hline\hline
\multicolumn{5}{|c|}{with oversampling, $w^+_i = w_i+2$} \\ \hline
 &
\multicolumn{2}{c|}{100 \%} &
\multicolumn{2}{c|}{25.2 \%} \\  \hline
1 	 & 0.506 & 0.716	 & 0.533 & 0.729 \\
2 	 & 0.223 & 0.472	 & 0.203 & 0.450 \\
4 	 & 0.069 & 0.258	 & 0.065 & 0.252 \\
8 	 & 0.025 & 0.143	 & 0.025 & 0.146 \\
12 & 0.010 & 0.083	 & 0.013 & 0.095 \\
16 & 0.006 & 0.058	 & 0.007 & 0.069 \\
\hline
\end{tabular}
\caption{Numerical tests for elasticity in heterogeneous media shown on the left of Figure \ref{fig:domain}. Randomized oversampling for GMsFEM.}
\label{tab:ov1}
\end{table}

\begin{table}
\centering
\begin{tabular}[hp]{|c|cc|cc|}
\hline
\multirow{2}{*}{$N_c$}  &
\multicolumn{2}{c|}{full snapshots} &
\multicolumn{2}{c|}{randomized snapshots} \\
 & $L_2$ & $H_1$ & $L_2$ & $H_1$ \\
 \hline \hline
\multicolumn{5}{|c|}{without oversampling, $w_i$} \\ \hline
 &
\multicolumn{2}{c|}{100 \%} &
\multicolumn{2}{c|}{39.2 \%} \\  \hline
1 	 & 0.692 & 0.822 & 0.858 & 0.908 \\
2 	 & 0.116 & 0.351 & 0.673 & 0.798 \\
4 	 & 0.039 & 0.207 & 0.510 & 0.678 \\
8 	 & 0.016 & 0.128 & 0.358 & 0.551 \\
12 & 0.007 & 0.078 & 0.204 & 0.418 \\
16 & 0.003 & 0.055 & 0.132 & 0.334 \\
\hline\hline
\multicolumn{5}{|c|}{with oversampling, $w^+_i = w_i+2$} \\ \hline
 &
\multicolumn{2}{c|}{100 \%} &
\multicolumn{2}{c|}{23.3 \%} \\  \hline
1 	 & 0.665 & 0.800	& 0.611 & 0.762 \\
2 	 & 0.115 & 0.355	& 0.131 & 0.383 \\
4   & 0.038 & 0.214	& 0.049 & 0.241 \\
8 	 & 0.015 & 0.135	& 0.020 & 0.157 \\
12 & 0.006 & 0.087	& 0.008 & 0.108 \\
16 & 0.006 & 0.089	& 0.008 & 0.102 \\
\hline
\end{tabular}
\caption{Numerical tests for elasticity in heterogeneous media shown on the right of Figure \ref{fig:domain}. Randomized oversampling for GMsFEM.}
\label{tab:ov2}
\end{table}

Comparing the results on the left columns (using full snapshot space) of Tables \ref{tab:ov1} and \ref{tab:ov2} with those on the right columns that correspond to the randomized snapshots, we observe that the randomized algorithm provides a nearly similar results, particularly with oversampling strategy. This is consistent with our previous observations for heterogeneous problems. Oversampling avoids the oscillations on the boundary conditions due to the randomness.
 The effect of oversampling strategy is much more obvious for the randomized snapshot space.

\section{Conclusion}
\label{sec:conclusion}

In this paper, we develop a generalized multiscale finite element framework
for problems in perforated domains. Our approach follows GMsFEM.
The main contributions of this paper are the development of snapshot space
and local spectral problems for solving problems in perforated domains
with multiple scales and no scale separation. Our approaches differ
from previously developed GMsFEM techniques for heterogeneous problems.
In particular, the snapshot vectors and local eigenvalue problems
need to be developed in a multiscale domain and take into account
the boundaries that may be disconnected. We show that
using GMsFEM framework, we can propose a unified framework for
solving problems in perforated domains. In the paper,
we discuss three applications:
(1) Laplace equation in perforated domain; (2) elasticity equation in
perforated domain; and (3) Stokes equations in perforated domain.
We present some preliminary numerical results
that show that one can efficiently solve these problems
on a coarse grid using fewer degrees of freedom. We also discuss the use
of randomized snapshots to reduce the offline computational cost
associated with computing the snapshot space.
In our future work, we plan to present analysis and design new more
efficient coarse spaces based on the analysis.

\section{Acknowledgement}
YE's work is partially supported by the U.S. Department of Energy Office of Science, Office of Advanced Scientific Computing Research, Applied Mathematics program under Award Number DE-FG02-13ER26165 and  the DoD Army ARO Project.

\bibliographystyle{siam}
\bibliography{references}
\end{document}